\documentclass[reqno]{amsart}
\usepackage{amssymb}
\usepackage{amscd}
\usepackage{amsmath}
\usepackage{amsthm}

\theoremstyle{plain}
  \newtheorem{theo}{Theorem}
  
  \newtheorem{prop}{Proposition}
  \newtheorem{lemm}{Lemma}
  \newtheorem*{quest*}{Question}  

\theoremstyle{remark}
  \newtheorem*{rema*}{Remark}
  \newtheorem{rema}{Remark}

\newcommand{\comment}[1]{}

\newcommand\ie{i.e.\ }

\newcommand\wrt{w.r.t.\ }

\renewcommand\o{\circ}

\newcommand\be{\beta}

\newcommand\de{\delta}

\newcommand\et{\eta}
\renewcommand\th{\theta}

\newcommand\ka{\kappa}
\newcommand\la{\lambda}
\newcommand\rh{\rho}

\newcommand\si{\sigma}
\newcommand\ta{\tau}
\newcommand\ph{\varphi}

\newcommand\om{\omega}
\newcommand\Ga{\Gamma}

\newcommand\La{\Lambda}
\newcommand\Si{\Sigma}
\newcommand\Om{\Omega}
\newcommand\Ph{\Phi}

\newcommand\ad{\text{\rm ad}}
\newcommand\Ad{\text{\rm Ad}}
\newcommand\pr{\text{\rm pr}}

\newcommand\Diff{\text{\rm Diff}}
\newcommand\oo{\infty}
\newcommand\id{\text{\rm id}}

\newcommand\ev{\text{\rm ev}}
\newcommand\flux{\text{\rm flux}}
\newcommand\per{\text{\rm per}}
\newcommand\Graph{\text{\rm Graph}}

\newcommand\sym{\text{\rm sym}}

\newcommand\g{\mathfrak g}

\newcommand\h{\mathfrak h}

\renewcommand\a{\mathfrak a}
\newcommand\X{\mathfrak X}

\newcommand\x{\times}

\newcommand\RR{\mathbb R}
\newcommand\ZZ{\mathbb Z}
\newcommand\TT{\mathbb T}

\begin{document}

\title[The path group construction]
{The path group construction\\ of Lie group extensions}

\author{Cornelia Vizman}

\address{Cornelia Vizman, 
         West University of Timisoara, Department of Mathematics, 
         Bd. V.Parvan 4, 300223-Timisoara, Romania.}

\email{vizman@math.uvt.ro}

\thanks{}



\begin{abstract}
We present an explicit realization of abelian extensions of 
infinite dimensional Lie groups using abelian extensions of path groups, 
by generalizing Mickelsson's approach to loop groups
and the approach of Losev-Moore-Nekrasov-Shatashvili
%\cite{LMNS98} 
to current groups. We apply our method to coupled cocycles on current Lie algebras and
to Lichnerowicz cocycles on the Lie algebra of divergence free vector fields.
\end{abstract}

\maketitle

%%%%%%%%%%%%%%%%%%%%%%%%%

\section{Introduction}

In infinite dimensions Lie's third theorem is not valid: there exist
Lie algebras which do not correspond to any Lie group.
In particular
given a connected infinite dimensional Lie group $G$,
not every abelian Lie algebra extension
$\a\to\hat\g\to\g$ of its Lie algebra $\g$
comes from a Lie group extension of $G$.
The obstructions determined in \cite {N04}
involve the fundamental and the second homotopy groups of $G$.
For instance if $G$ is simply connected, the integrability condition 
for the Lie algebra extension $\hat\g$ 
described by the Lie algebra 2-cocycle $\om$ on $\g$ with values in 
the $\g$-module $\a$,
reduces to the discreteness of the period group $\Pi_\om\subset\a$
of the cocycle $\om$.
Under this assumption, for any discrete subgroup $\Ga$ of the space $\a^G$
of $G$-invariant elements of $\a$, containing the period group,
%with $\Pi_\om\subseteq\Ga$,
there exists a corresponding abelian Lie group extension $A\to\hat G\to G$ of $G$ by 
the $G$-module $A=\a/\Ga$. 

%Important for field theory is the 
Much studied is the central extension of the loop group
$G=C^\oo(S^1,H)$ of a simple Lie group $H$.
With a suitable multiple $\ka$ of the Killing form 
of the Lie algebra $\h$ of $H$,
a Lie algebra 2-cocycle on the loop algebra 
$\g=C^\oo(S^1,\h)$ is 
$$
\om(X,Y)=2\int_{S^1}\ka(X,dY). 
$$
When the simple Lie group $H$ is simply connected 
then $G$ is also simply connected, and the period group is $\Pi_\om=\ZZ$,
so $\om$ is integrable.
Explicit constructions of the corresponding central extension
$\TT\to\hat G\to G$
of the loop group $G$ by the circle $\TT=\RR/\ZZ$
can be found in Chapter 4 of \cite{PS86}, in \cite{M87}, \cite{M88}, and in Chapter 4 of \cite{M89}. 

The approach in the book of Pressley-Segal \cite{PS86} is very general:
one considers a simply connected 
prequantizable manifold $(M,\Om)$ (\ie $\Om$ is a closed integral 2-form on $M$
and there is a principal circle bundle $P$ over $M$ 
with a principal connection having curvature $\Om$),
together with a $G$-action preserving $\Om$
and one pulls back Kostant's prequantization central extension \cite{K70}
to $G$ by this action.
More precisely the resulting extension $\hat G$ 
is the group of all fiber preserving diffeomorphisms
of $P$ which preserve the connection 1-form and cover an element of $G$.
A concrete description of the central extension of $G$ obtained in this way,
using paths in $M$, is given there.
Central Lie group extensions associated to Hamiltonian
actions on a prequantizable manifold $(M,\Om)$  
are considered in \cite{NV03}.
%, the action being called Hamiltonian if the insertion of any fundamental vector field in $\Om$ gives an exact 1-form. A flux homomorphism in this setting is also needed
A generalization of this approach to abelian extensions
is the subject of \cite{V06}.

A second approach \cite {M87} and \cite{M89} is due to Mickelsson: 
the central extension of the loop group
is explicitly realized as a quotient group of a central extension of the group
of currents on the 2-disk $D$. The last central extension is given by the group 2-cocycle 
$c(f,g)=\int_D\de^lf\wedge_\ka\de^rg$, where $\de^l$ and $\de^r$
denote the left and right logarithmic derivative. 
This construction is generalized to central extensions of current groups
on Riemann surfaces in \cite{FK96}. The construction of central extensions
of current groups on arbitrary compact manifolds \cite{LMNS98}
is due to Losev-Moore-Nekrasov-Shatashvili
explicitly using the path group of the current group.

In this paper we generalize this last approach,
obtaining a path group method for the construction of
abelian Lie group extensions.
We consider a connected Lie group $G$
and the exact sequence of Lie groups
$\Om_0G\to PG\to\tilde G$ where $\tilde G$ is the universal covering group,
$PG$ the path group and $\Om_0G$ its subgroup of null-homotopic loops.
We also consider a 2-cocycle $\om$ on $\g$ with values in the $\g$-module
$\a$, having a discrete period group,
and the discrete subgroup $\Ga\supseteq\Pi_\om$ of $\a^G$ with $A=\a/\Ga$.
An abelian Lie group extension of $\tilde G$ by $A$ integrating $\a\rtimes_\om\g$
%the abelian Lie algebra extension of $\g$ defined by $\om$, 
is obtained as the quotient of an 
abelian Lie group extension $A\rtimes_cPG$,
% defined by a group cocycle $c:PG\x PG\to A$,
by the graph of a map $\la:\Om_0G\to A$.
The group cocycle $c:PG\x PG\to A$, as well as the map $\la$
are explicitly given by formulas (Theorem \ref{main}).
%(\ref{11}) and (\ref{12}).
%This is the content of Section 2 and 3.
A geometric construction of abelian Lie group extensions using
the path group is presented in \cite{H07}.
We show that the two group extensions are isomorphic.
%Pressley-Segal and to that of Hekmati.

Our method works well in concrete settings, in spite of the heaviness
of the formula for $c$.
In Section 4 we specialize to loop groups and current groups.
Here a symmetrization procedure applied to $c$ and $\la$
simplify considerably the result, thus recovering
the constructions in \cite{M89} and \cite{LMNS98}.
In Section 5 we treat the coupled cocycle on current Lie algebras 
defined by Neeb in \cite{N07}.
%from a work in progress by Neeb.
%\cite{N07}, \cite{NW07}.
Finally we use this construction in Section 6 to explicitly realize central
extensions of the group of volume preserving diffeomorphisms
integrating Lichnerowicz cocycles.

{\bf Acknowledgements:} I am grateful to Karl-Hermann Neeb for showing 
me a variant of Theorem \ref{main} belonging to an extended version of \cite{N02},
for the coupled cocycle and for useful comments.

%%%%%%%%%%%%%%%%%%%%%%%%

%%%%%%%%%%%%%%%%%%%%%%%%%

\section{The path group}\label{2}

Let $G$ be a Lie group, $\g$ its Lie algebra, 
$\a$ a smooth $G$-module (\ie the action map $\rh_\a:G\x\a\to\a$ is smooth) 
and $\a^G$ the subspace of $G$-invariant elements.
The continuous $\a$-valued Lie algebra 2-cocycle $\om$ on $\g$ defines a closed
equivariant $\a$-valued 2-form $\om^{eq}$ on $G$.
This means $\rh_\a(g)\o\om^{eq}=\la_g^*\om^{eq}$, $\la_g$ 
denoting the left translation by $g$ in $G$. 
The period homomorphism of the 2-cocycle $\om$ is by definition 
\begin{equation}\label{period}
\per_\om:\pi_2(G)\to\a^G,\quad
\per_\om([\ta])=\int_{S^2}\ta^*\om^{eq}
\end{equation} 
for $\ta:S^2\to G$ a piecewise smooth representative of a free homotopy class
(a spherical 2-cycle).
The image $\Pi_\om\subset\a^G$ of $\per_\om$ is called the period group of $\om$.

We present van Est's method \cite{E54} for obtaining group 2-cocycles by integrating
the closed 2-form $\om^{eq}$ over suitable triangles.
Let $G$ be a simply connected Lie group and $\{a_g\}_{g\in G}$
a given family of smooth paths, $a_g:I=[0,1]\to G$ from $e$ to $g$.
For each $f,g\in G$ let $\Si_{f,g}$ be a piecewise smooth 2-simplex in $G$ 
with boundary $fa_g-a_{fg}+a_f$ and let $C$ be the map 
\begin{equation}\label{c}
C:G\x G\to\a,\quad C(f,g)=\int_{\Si_{f,g}}\om^{eq}.
\end{equation}
Because $\th_{f,g,h}=f\Si_{g,h}-\Si_{fg,h}+\Si_{f,gh}-\Si_{f,g}$
is a spherical 2-cycle in $G$ for any $f,g,h\in G$, and because $\om^{eq}$ is equivariant, 
the map $C$ satisfies the relation
\begin{equation}
%\label{cocycle}
fC(g,h)-C(fg,h)+C(f,gh)-C(f,g)
%&=\int_{\Si_{g,h}}f\om^{eq}-\int_{\Si_{fg,h}}\om^{eq}+\int_{\Si_{f,gh}}\om^{eq}
%-\int_{\Si_{f,g}}\om^{eq}\\
=\int_{\th_{f,g,h}}\om^{eq}\in\Pi_\om.
\end{equation}

If the period group $\Pi_\om$ of the Lie algebra 2-cocycle $\om$ is discrete
and $\Ga\supseteq\Pi_\om$ is any discrete subgroup of the space $\a^G$ 
of $G$-invariant elements, 
we denote the abelian Lie group and smooth $G$-module $\a/\Ga$ 
by $A$ and the quotient map by $\exp:\a\to A$.
Then $c=\exp\o C:G\x G\to A$ is a group 2-cocycle 
%by (\ref{cocycle}).
independent of the choice of the 2-cycles $\Si_{f,g}$.
When the paths $a_g$ are carefully chosen \cite{N04}, then $c$ is smooth in an identity
neighborhood and $\om$ is the associated Lie algebra cocycle,
\ie $d^2c_{(e,e)}(X,Y)-d^2c_{(e,e)}(Y,X)=\om(X,Y)$ for all $X,Y\in\g$.

\begin{rema}\label{vanest}
Given a smoothly contractible Lie group $G$, 
each smooth retraction $h:I\x G\to G$ of $G$ to $\{e\}$ provides 
a family of smooth paths $\{a_g\}$ from $e$ to $g$ by 
%\begin{equation}\label{homotopy}
$a_g(s)=h(s,g), s\in I$.
%\end{equation} 
Any continuous Lie algebra 2-cocycle $\om\in Z ^2(\g,\a)$ 
is integrable to a smooth $\a$-valued van Est cocycle $C$ on $G$ 
given by (\ref{c}), depending only on the retraction $h$.
One can choose for instance $\Si_{f,g}(s,t)=h(s,fh(t,g))$, $s,t\in I$.
\end{rema}

%%%%%%%

Let $G$ be a connected Lie group with Lie algebra $\g$.
The group of smooth paths in $G$ starting at the identity, 
$$
PG=\{g\in C^\oo(I,G)|g(0)=e\},
$$ 
called the path group, is a 
smoothly contractible Lie group
with Lie algebra 
$$
P\g=\{X\in C^\oo(I,\g)|X(0)=0\},
$$ 
the path Lie algebra.

Each $\a$-valued Lie algebra 2-cocycle $\om$ on $\g$ defines an 
$\a$-valued Lie algebra 2-cocycle 
$P\om=\ev_1^*\om:(X,Y)\mapsto\om(X(1),Y(1))$ on the path Lie algebra
$P\g$, the evaluation map $\ev_1:P\g\to\g$ being a Lie algebra homomorphism.
Via the group homomorphism $\ev_1:PG\to G$, the $G$-module $\a$
becomes a $PG$-module.
With Remark \ref{vanest} we can integrate $P\om$ to 
a smooth group 2-cocycle $C$ on the contractible group $PG$.
One can write down an explicit formula for this cocycle
using the left logarithmic derivative $\de^lf\in P\g$ for $f\in PG$.

\begin{prop}\label{central}
A smooth group 2-cocycle on the path group $PG$ integrating the Lie algebra 
cocycle $P\om$ is
\begin{equation}\label{formula}
C(f,g)=\int_0^1\Big(\int_0^s\rh_\a(f(s)g(t))\om(\Ad(g(t)^{-1})\de^lf(s),\de^lg(t))dt\Big)ds.
\end{equation}
\end{prop}

\begin{proof}
A smooth retraction of $PG$ to $e$ is obtained by reparametrisation
of paths, namely $h:I\x PG\to PG$,  
$h(s,g)(t)=g(st)$ for $g\in PG$ and $s,t\in I$.
The system of paths $\{a_g\}_{g\in PG}$ in $PG$ defined by $h$ 
is $a_g(s)(t)=g(st)$. It
has the property $a_fa_g=a_{fg}$ for all $f,g\in PG$, so we choose a particular
2-simplex $\Si_{f,g}$ in $PG$ with boundary $fa_g-a_{fg}+a_f$, 
namely $\Si_{f,g}(s,t)=a_f(s)a_g(st)$ for $(s,t)\in I\x I$.
%If the equivariant 2-form on $G$ defined by $\om$ is $\om^{eq}$, 
The equivariant 2-form on $PG$ defined by $P\om$ is
$\ev_1^*\om^{eq}$ and the van Est cocycle on $PG$ integrating $P\om$ is
$C(f,g)=\int_{\Si_{f,g}}\ev_1^*\om^{eq}$.

The 2-simplex $\si_{f,g}=\ev_1\o\Si_{f,g}$ on $G$ is given by
$\si_{f,g}(s,t)=f(s)g(st)$ for $(s,t)\in I\x I$, hence
\begin{align*}
C(f,g)&=\int_{\si_{f,g}}\om^{eq}
%=\int_0^1\int_0^1\om^{eq}(\partial_s\si_{f,g}(s,t),\partial_{t}\si_{f,g}(s,t))dtds\\
=\int_0^1\int_0^1\om^{eq}(\dot f(s)g(st)+tf(s)\dot g(st),sf(s)\dot g(st))dtds\\
&=\int_0^1\int_0^1\rh_\a(f(s)g(st))\om(\Ad(g(st)^{-1})\de^lf(s)
+t\de^lg(st),s\de^lg(st))dtds\\
&=\int_0^1\Big(\int_0^s\rh_\a(f(s)g(t))\om(\Ad(g(t)^{-1})\de^lf(s),\de^lg(t))dt\Big)ds,
\end{align*}
using the $G$-equivariance of $\om^{eq}$ at step 3
and a change of variable at step 4.
\end{proof}

\begin{rema}\label{prime}
Another group cocycle on $PG$ integrating the Lie algebra 2-cocycle
$P\om$ is $C'(f,g)=-C(g^{-1},f^{-1})$.
It can be seen as the van Est cocycle associated to the system of paths 
$\{b_g\}_{g\in PG}$, $b_g(s)=ga_g(1-s)^{-1}$.
From (\ref{formula}) we obtain in particular
\begin{equation}\label{cprime}
C'(f,g)=-\int_0^1\Big(\int_0^s\rh_\a(g(t)^{-1}f(s)^{-1})
\om(\Ad(f(t))\de^rg(s),\de^rf(t))dt\Big)ds.
\end{equation}

By the symmetrization procedure
\begin{equation}\label{csym}
C_{\sym}(f,g)=\frac12(C(f,g)+C'(f,g))=\frac12(C(f,g)-C(g^{-1},f^{-1}))
\end{equation}
we get a new group cocycle $C_{\sym}$ on $PG$ integrating $P\om$
and having the property $C_{\sym}(g^{-1},f^{-1})=-C_{\sym}(f,g)$.
%For a general group cocycle $c$ on $G$ with values in the abelian group $A$
%and integrating $\om$, we can build the group cocycle 
%$c_s(f,g)=c(f,g)c(g^{-1},f^{-1})^{-1}$ integrating $2\om$.
\end{rema}

\begin{rema}
When $\a$ is a trivial $G$-module and $\om^l$ the left invariant 2-form on $G$
defined by $\om\in Z^2 (\g,\a)$,
then a formula for the van Est group 2-cocycle $C$ on $PG$ integrating $P\om$ is
\begin{equation}\label{formula2}
C(f,g)=\int_{\si_{f,g}}\om^l
=\int_0^1\Big(\int_0^s\om(\Ad(g(t)^{-1})\de^lf(s),\de^lg(t))dt\Big)ds.
\end{equation}
\end{rema}

\begin{rema}
Given a continuous Lie algebra $p$-cocycle $\om$ on $\g$ with values
in the smooth $G$-module $\a$, a group $p$-cocycle $C$ on $PG$ integrating $P\om$
is
\begin{align*}
C(g_1,\dots,g_p)
%&=\int_{\si_{g_1,\dots,g_p}}\om^{eq}\\
&=\int_0^1\Big(\int_0^{t_1}\dots\Big(\int_0^{t_{p-1}}\rh_\a(g_1(t_1)\dots g_p(t_p))
\cdot\om(\Ad(g_2(t_2)g_3(t_3)\dots \\
&g_p(t_p))^{-1}\de^lg_1(t_1),
\Ad(g_3(t_3)\dots 
g_p(t_p))^{-1}\de^lg_2(t_2),\dots,\\
&\Ad(g_p(t_p))^{-1}\de^lg_{p-1}(t_{p-1}),\de^lg_p(t_p))dt_p\Big)\dots dt_2\Big)dt_1,
\end{align*}
obtained by integrating the closed equivariant $p$-form $\om^{eq}$ on $G$
over the $p$-simplex $\si_{g_1,\dots,g_p}:(t_1,\dots,t_p)\in I^p\mapsto
g_1(t_1)g_2(t_1t_2)\dots g_p(t_1t_2\dots t_p)$.
\end{rema}

%%%%%%%%%%%%%%%%%%%%%%%%%%%%%%%%%%%%%%%%%%%%%%%%%%%

\section{Construction of abelian Lie group extensions via path groups}

Let $\om$ be an $\a$-valued Lie algebra 2-cocycle on the Lie algebra $\g$
of the connected Lie group $G$. 
We assume that its period group $\Pi_\om\subset\a^G$ is discrete and 
$\Ga\supseteq\Pi_\om$ is a discrete subgroup of $\a^G$,
and we denote by $\exp:\a\to A=\a/\Ga$ the quotient map.
The smooth $G$-action $\rh_\a$ on $\a$
descends to a smooth $G$-action $\rh_A$ on $A$,
because $\Ga\subset\a^G$ .

In this section we explicitly realize an abelian extension
of the universal covering group $\tilde G$ of $G$ by the abelian Lie group
$A$. It is done by factorising an abelian extension 
of the path group $PG$ (defined with a van Est cocycle) 
by the graph of a smooth map.

The subgroup $\Om G\subset PG$ of loops based at $e$ has a subgroup
$\Om_0 G\subseteq \Om G$ of null-homotopic loops based at $e$. 
Both have as Lie algebra
the Lie algebra $\Om\g$ of loops in $\g$ based at $0$,
moreover $\Om_0G$ is the identity component of $\Om G$.
The following two sequences of Lie groups are exact:
$$
\CD
\Om G & @>>> & PG & @>>> & G &
\\
@AAA & & @| & & @AAA &
\\
\Om_0G & @>>> & PG & @>>> & \tilde G
\endCD
$$
\comment
{
\begin{gather*}
\Om G\to PG\to G\\
\Om_0 G\to PG\to\tilde G
\end{gather*} 
}
and have the same exact sequence of Lie algebras
$\Om\g\to P\g\to\g$.

Let $C:PG\x PG\to\a$ be the van Est 2-cocycle (\ref{formula}) integrating the Lie algebra 2-cocycle $P\om$ on $P\g$. 
The cocycle $c=\exp\o C:PG\x PG\to A$ is given by
\begin{equation}\label{11}
c(f,g)
=\exp\int_{\si_{f,g}}\om^{eq}
=\exp\int_0^1\Big(\int_0^s\rh_\a(f(s)g(t))\om(\Ad(g(t)^{-1})\de^lf(s),\de^lg(t))dt\Big)ds,
\end{equation}
where $\si_{f,g}(s,t)=f(s)g(st)$, $s,t\in I$ for $f,g\in PG$,
The Lie algebra
cocycle $P\om$ vanishes on $\Om\g$.
The next proposition will show that the restriction of the group cocycle $c$ to $\Om_0G$ is a coboundary (both $\a$ and $A$ are trivial $\Om_0G$-modules).

We consider the smooth map
\begin{equation}\label{big}
\La:\widetilde{\Om_0G}\to\a,\quad\La([\bar g])=-\int_{\bar g}\om^{eq},
\end{equation}
viewing the path $\bar g$ from $e$ to $g$ in $\Om_0G$ as a map $\bar g:I\x I\to G$.
It is well defined since the integral of the closed 2-form $\om^{eq}$ 
over two homotopic paths 
in $\Om_0G$ (hence homotopic maps from $I\x I$ to $G$) is the same.

\begin{rema}\label{periods}
Identifying $\pi_2(G)$ with $\pi_1(\Om_0G)\subset\widetilde{\Om_0G}$,
the restriction of $\La$ to $\pi_1(\Om_0G)$ equals $-\per_\om$,
the opposite of the period map (\ref{period}). Indeed, a loop at $e$
of loops in $G$ determines a spherical 2-cycle in $G$.
In particular the map $\La$ is well defined on $\Om_0G$
when considered modulo $\Pi_\om$,
hence it descends to a well defined map
%with values in $\a/\Pi_\om$: 
\begin{equation}\label{lambda}
\la:\Om_0G\to A,\quad\la(g)=\Big(\exp\int_{\bar g}\om^{eq}\Big)^{-1},
\end{equation}
\ie $\la(g)$ does not depend on the chosen path $\bar g$ in $\Om_0G$ from $e$ to $g$.
\end{rema}

\begin{prop}\label{coboundary}
The identity
\begin{equation}\label{delambda}
c(f,g)=\la(fg)\la(f)^{-1}\la(g)^{-1}
\end{equation} 
holds for all $f,g\in\Om_0G$. In particular $c$ restricted to $\Om_0G$
is the coboundary of $\la^{-1}$.
\end{prop}

\begin{proof}
Let $\bar f$ and $\bar g$ be paths in $\Om_0G$ with $f=\bar f(1)$ and $g=\bar g(1)$.
The 2-chains 
\begin{equation*}
\si_{f,g}:I\x I\to G,
\quad\si_{f,g}(s,t)=f(s)g(st)
\end{equation*}
and $\bar f+\bar g-\bar f\bar g$
have the same boundary $f+g-fg$, so they determine a spherical 2-cycle 
$\ta_{f,g}=\si_{f,g}-\bar f-\bar g+\bar f\bar g$ in $G$.
Hence the van Est cocycle $C(f,g)=\int_{\si_{f,g}}\om^{eq}$ on $PG$ satisfies
\begin{equation}\label{cla}
C(f,g)+\La([\bar f])+\La([\bar g])-\La([\bar f\bar g])=\int_{\ta_{f,g}}\om^{eq}\in\Pi_\om
\subseteq\Ga.
\end{equation}
% (from the proof of Proposition \ref{central}).
The projection of this relation to $A$ gives the requested identity.
\end{proof}

\begin{rema}\label{label}
The stronger relation
\begin{equation}\label{cla1}
C(f,g)+\La([\bar f])+\La([\bar g])-\La([\bar f\bar g])=0,\quad f,g\in\Om_0G
\end{equation}
also holds, because there always exists a bordism from the 2-cycle $\ta_{f,g}$ to $e$
given by $(r,t,s)\in I\x I\x I\mapsto \bar f(r,s)\bar g(r,st)
-\bar f(rt,s)-\bar g(rs,t)+\bar f(rt,s)\bar g(rt,s)\in G$, for
$\bar f$ and $\bar g$ paths in $\Om_0G$.
Still the $\a$-valued van Est cocycle $C$ restricted to $\Om_0G$ is not a coboundary
in general:
$\La$ does not descend to a well defined $\a$-valued map on $\Om_0G$. 
Anyway, following \cite{FK96}, if (\ref{cla1}) is satisfied
we say that the map $\La$ resolves the 2-cocycle $C$.
\end{rema}

\begin{rema}
The cocycles $C'$ and $C_{\sym}$ defined in Remark \ref{prime}
also posess resolving maps.
The map 
\begin{equation}\label{bigprime}
\La'([\bar g])=-\La([\bar g^{-1}])=\int_{\bar g^{-1}}\om^{eq}
\end{equation}
resolves the cocycle $C'$ 
and $\La_{\sym}=\frac12(\La+\La')$ resolves the cocycle $C_{\sym}$.
In particular we have that $\La_{\sym}([\bar g]^{-1})=-\La_{\sym}([\bar g])$. 
\end{rema}

\begin{lemm}\label{lemma}
Let $H$ be a normal split Lie subgroup of the Lie group $G$
and $A$ a smooth $G$-module, trivial as an $H$-module.
Let $c$ be an $A$-valued group 2-cocycle on the   group $G$
whose restriction to $H$ is the coboundary of $\la^{-1}$ for 
a given smooth map $\la:H\to A$.
When one of the following two equivalent conditions:
\begin{enumerate}
\item $c(g,h)c(ghg^{-1},g)^{-1}=(\rh_A(g)\la(h)^{-1})\la(ghg^{-1})$
for all $g\in G$ and $h\in H$
\item the graph of $\la$ is a normal subgroup of $A\rtimes_c G$
\end{enumerate}
is satisfied, then the quotient group $(A\rtimes_c G)/{\Graph(\la)}$ 
is an abelian Lie group extension of $G/H$ by $A$.
\end{lemm}

\begin{proof}
The graph of $\la$ coincides with the image of
the map $\ph:h\in H\mapsto (\la(h),h)\in A\rtimes_c G$.
Because $c$ is the coboundary of $\la$:
\begin{equation}\label{h1h2}
c(h_1,h_2)=\la(h_1h_2)\la(h_1)^{-1}\la(h_2)^{-1}\text{ for all }h_1,h_2\in H,
\end{equation}
$\ph$ is a   group homomorphism. This follows from
$$
\ph(h_1)\ph(h_2)
=(\la(h_1)\la(h_2)c(h_1,h_2),h_1h_2)
=\ph(h_1h_2).
$$
Hence the graph of $\la$ is a subgroup of $A\rtimes_c G$.

%When both $g,h\in H$, then relation (i) follows from (\ref{h1h2}).
Let $g\in G$ and $h\in H$. 
The conjugate in $A\rtimes_c G$ of the element $(\la(h),h)\in\Graph(\la)$ is
(see for instance Lemma 2.1 in \cite{N04})
$$
(a,g)(\la(h),h)(a,g)^{-1}=((\rh_A(g)\la(h))c(g,h)c(ghg^{-1},g)^{-1},ghg^{-1}).
$$
It belongs to the graph of $\la$, \ie it equals $(\la(ghg^{-1}),ghg^{-1})$, 
if and only if the identity (i) holds.

The kernel of the projection homomorphism $(a,g)\in (A\rtimes_cG)/
{\Graph(\la)}\mapsto gH\in G/H$ is isomorphic to $A$,
hence $(A\rtimes_cG)/{\Graph(\la)}$ is an abelian   group extension
of $G/H$ by $A$, for the natural $G/H$-module structure on $A$. 
%($H$ acting trivially on $A$).
It is a Lie group since $\Graph(\la)$ is a split Lie subgroup of $A\rtimes_cG$,
$H$ being a split Lie subgroup of $G$.
\end{proof}

\begin{theo}\label{main}
Let $G$ be a connected Lie group, $\a$ a smooth $G$-module and $\om\in Z^2 (\g,\a)$
a continuous Lie algebra 2-cocycle with discrete period group $\Pi_\om$.
Let $\Ga\supseteq\Pi_\om$ be a discrete subgroup of $\a^G$, $\exp:\a\to A=\a/\Ga$
the quotient map and $c$ the cocycle defined in (\ref{11}).

Then the graph of the smooth map $\la$ defined in (\ref{lambda})
is a normal subgroup of $A\rtimes_cPG$ and
the quotient group $(A\rtimes_cPG)/\Graph(\la)$
is an abelian Lie group extension of the universal covering group $\tilde G$ 
by $A$, integrating the Lie algebra extension $\a\rtimes_\om\g$.
\end{theo}

\begin{proof}
To apply Lemma \ref{lemma} to the normal split Lie subgroup
$\Om_0G$ of $PG$,
we verify relation (i) for the $A$-valued 2-cocycle $c$ on $PG$
and the map $\la:\Om_0G\to A$. 

The boundary of the 2-chain $\si_{f,g}$ for $f,g\in PG$
is $f(1)g-fg+f$ and the boundary of the 2-chain $\bar f$ for $f\in\Om_0G$
is $f$.
Let $g\in PG$ and $h\in\Om_0G$ with $\bar h$
a path in $\Om_0G$ from $e$ to $h$. 
The 2-chains 
$\si_{ghg^{-1},g}-\si_{g,h}$ and $g\bar hg^{-1}-g(1)\bar h$ in $G$ have
the same boundary, namely $ghg^{-1}-g(1)h$.
Integrating $\om^{eq}$ over the spherical 2-cycle $\nu_{g,h}=\si_{ghg^{-1},g}-\si_{g,h}
-g\bar hg^{-1}+g(1)\bar h$ and using the $G$-equivariance of $\om^{eq}$ we obtain
$$
C(ghg^{-1},g)-C(g,h)+\La(g\bar hg^{-1})-\rh_A(g(1))\La(\bar h)
=\int_{\nu_{g,h}}\om^{eq}\in\Pi_\om\subseteq\Ga.
$$
The projection of this identity to $A$ gives (i),
showing that the graph of $\la$ is a normal subgroup of $A\rtimes_cPG$.

The abelian Lie group extension $A\rtimes_cPG$ of $PG$ integrates $P\om=\ev_1^*\om$,
hence the quotient group $(A\rtimes_cPG)/\Graph(\la)$
is an abelian Lie group extension of 
the universal covering group $\tilde G=PG/\Om_0G$ integrating $\om$. 
\end{proof}

The rows and the last column of the
following diagram are exact sequences of Lie groups
$$
\CD
A & @>>> & (A\rtimes_cPG)/\Graph(\la) & @>>>  & \tilde G &
\\
@AAA &                & @AAA &  & @AAA &
\\
A & @>>> & A\rtimes_cPG & @>>> & PG &
\\
@AAA &                & @AAA &  & @AAA &
\\
A  & @>>> &   A\times_c\Om_0G  & @>>>  &  \Om_0G.
\endCD
$$

\begin{rema}
There is a geometric construction of an abelian extension
of $\tilde G$ using the path group $PG$ presented in \cite{H07}.  
%and generalizing the construction in \cite{PS86}. 
One considers the following equivalence relation on $A\rtimes_cPG$,
where $c$ is given by (\ref{11}):
\begin{equation}\label{simm}
(a_1,g_1)\sim(a_2,g_2)\Longleftrightarrow g_1(1)=g_2(1), g_1-g_2=\partial\si,
a_2=a_1\exp\int_\si\om^{eq},
\end{equation}
the second condition on the right meaning that $\si$ is any 2-chain in $G$ having as boundary the loop $g_1-g_2$.
Then $ (A\rtimes_cPG)/\sim$ is an abelian extension of $\tilde G$ integrating $\om$.
We show it is isomorphic to the abelian extension in Theorem \ref{main}.

%In particular $(a,g)\sim(1,e)$ if and only if
%$g\in\Om_0G$ is a contractible loop and $a=(\exp\int_{\bar g}\om^{eq})^{-1}$
%for a path $\bar g$ in $\Om_0G$ from $e$ to $g$, 
%\ie $a=\la(g)$ by (\ref{12}).
A pair $(a,g)$ is equivalent to the identity element if and only if 
$(a,g)$ belongs to $\Graph(\la)$. 
Moreover, two pairs $(a_1,g_1)$ and $(a_2,g_2)$
are equivalent if and only if the composition $(a_1,g_1)^{-1}(a_2,g_2)$
taken in $A\rtimes_cPG$ belongs to $\Graph(\la)$.
Indeed, for $h\in\Om_0G$ and $\bar h$ a path from $e$ to $h$ in $\Om_0G$,
\begin{align*}
(a_1,g_1)(\la(h),h)
%&=\Big(a_1\rh_A(g_1(1))
%\Big(\exp\int_{\bar h}\om^{eq}\Big)^{-1}c(g_1,h),g_1h\Big)\\
&=\Big(a_1\Big(\exp\int_{g_1(1)\bar h}\om^{eq}\Big)^{-1}
\exp\int_{\si_{g_1,h}}\om^{eq},g_1h\Big)\\
&=\Big(a_1\exp\Big(\int_{\si_{g_1,h}-g_1(1)\bar h}\om^{eq}\Big),g_1h\Big),
\end{align*}
so $(a_1,g_1)(\la(h),h)=(a_2,g_2)$ if and only if $g_2=g_1h$ and $a_2=a_1\exp\int_\si\om^{eq}$,
where $\si$ is any 2-cycle in $G$ such that
$\partial\si=\partial(\si_{g_1,h}-g_1(1)\bar h)=g_1-g_1h=g_1-g_2$.
Hence the abelian extension $(A\rtimes_cPG)/\Graph(\la)$ from
Theorem \ref{main} and the abelian extension 
$(A\rtimes_cPG)/_\sim$ are isomorphic.
\end{rema}

%%%%%%%%%%%%%%%%%%%%%%%%%%
%%%%%%%%%%%%%%%%%%%%%

\section{Current groups}

Let $M$ be a compact manifold and $H$ a finite dimensional
connected Lie group with Lie algebra $\h$.
The current group $C^\oo(M,H)$ with pointwise multiplication
is a Lie group with Lie algebra the current algebra $\g=C^\oo(M,\h)$ 
(as in \cite{KM97} Section 42).
By $G$ we denote the identity component $C^\oo(M,H)_0$ of the current group.

We consider an invariant symmetric bilinear form
$$
\ka:\h\x\h\to V.
$$
Defining $\a=\Om^1(M,V)/dC^\oo(M,V)$,
there is a continuous Lie algebra 2-cocycle on the current algebra,
\begin{equation}\label{mn}
\om:\g\x\g\to\a,\quad\om(X,Y)=[\ka(X,dY)-\ka(Y,dX)]=2[\ka(X,dY)].
\end{equation}  
In the loop group case $M=S^1$ the space $\a$ can be identified with $V$, 
so the cocycle on the loop algebra $\g=C^\oo(S^1,\h)$ can be taken as
\begin{equation}\label{km}
\om:\g\x\g\to V,\quad\om(X,Y)=2\int_{S^1}\ka(X,dY).
\end{equation}

\begin{rema}\label{discrete}
When $H$ is simply connected and simple, the loop group 
$G$ is simply connected. If $\ka:\h\x\h\to\RR$ is the suitably normalized
Killing form, $V=\RR$ and the period group is $\Pi_\om=\ZZ\subset\RR$.
In Chapter 4 of \cite{M89} is presented
the construction of the central extension of $G$ by the circle
$\TT=\RR/\ZZ$, integrating $\om$.
It is explicitly realized as a quotient group of the central extension of the group
of currents on the 2-disk $D$ given by the group 2-cocycle 
$c(f,g)=\int_D\de^lf\wedge_\ka\de^rg$.

The generalization of this result to current groups can be found in \cite{LMNS98}.
The bilinear form $\ka$ is a multiple of the Killing form of the simple Lie algebra $\h$, and
the space
%There $\h$ is a simple Lie algebra and $\ka$ is a multiple of the Killing form on $\h$,
$\a$ is $\Om^1(M)/dC^\oo(M)$.
When $\ka$ is suitably normalized, then the period group 
is the discrete subgroup $\Pi_\om=Z^1_{\ZZ}(M)/dC^\oo(M)$ of $\a$ consisting
of all cohomology classes with integral periods.
%and it is contained in the subspace $H^1_{dR}(M)$ of $\a$. 
The central extension of the universal covering group $\tilde G$ 
by $A=\a/\Pi_\om=\Om^1(M)/Z^1_{\ZZ}(M)$
integrating $\om$ is constructed in \cite{LMNS98}
as a quotient group of the central extension of the path group of $G$
given by the group cocycle $c(f,g)=[\int_I\de^lf\wedge_\ka\de^rg]$, $f,g:I\x M\to H$.
Adapting this construction one gets also the central extension
of the gauge group of automorphisms of a nontrivial vector bundle
\cite{LMNS98}.
\end{rema}

\begin{rema}
The period group of $\om$ is not always discrete; an example with non-discrete period group 
can be found in Remark II.10 \cite{MN03}. The reduction theorem I.6 \cite{MN03}
shows that given a $V$-valued invariant symmetric form $\ka$,
the discreteness of the period group for $M=S^1$ implies the discreteness of the period group
for any compact manifold $M$.
\end{rema}

\begin{rema}
Considering the natural action of $\h$ on the symmetric power $S^2(\h)$
induced by the adjoint action, the quotient space $V(\h)=S^2(\h)/(\h\cdot S^2(\h))$
comes with a universal invariant symmetric bilinear form
$\ka_\h:\h\x\h\to V(\h)$.
For any invariant symmetric bilinear form $\ka$ on $\h$ with values in $V$,
there is a unique linear map $\ph:V(\h)\to V$ such that $\ka=\ph\o\ka_\h$.
The continuous Lie algebra 2-cocycle
$\om_\h(X,Y)=2[\ka_\h(X,dY)]$ on the current algebra
with values in $\Om^1(M,V(\h))/dC^\oo(M,V(\h))$
has a discrete period group $\Pi_\om$ contained in the subspace
$H^1_{dR}(M,V(\h))$.
\cite{MN03}.    
\end{rema}

Assuming the period group $\Pi_\om\subset\a$ of the continuous Lie  algebra
2-cocycle (\ref{mn}) is discrete,
we apply Theorem \ref{main} to the identity component of the current group
to explicitly realize a central extension of its universal covering group integrating $\om$.

The path group of $G=C^\oo(M,H)_0$ and the path Lie algebra of $\g=C^\oo(M,\h)$ are
(by \cite{KM97} Section 42):
\begin{gather*}
PG=\{g\in C^\oo(I\x M,H)|g(0,x)=e,\forall x\in M\}=PC^\oo(M,H)\\
P\g=\{X\in C^\oo(I\x M,\h)|X(0,x)=0,\forall x\in M\},
\end{gather*}
because the path group of a Lie group coincides with the path group of
its identity component. 
Denoting by $d_xX$ the exterior differential of $X$ on $M$
and by $\de^l_xg$ the logarithmic derivative of $g$ on $M$
(viewing $t\in I$ as a parameter),
we get 
$$
dX\wedge_\ka dY=dt\wedge(\ka(\frac{d}{dt}X,d_xY)
-\ka(d_xX,\frac{d}{dt}Y))+d_xX\wedge_\ka d_xY.
$$
An expression for the pullback 2-cocycle $P\om=\ev_1^*\om$ on $P\g$ is in this case
\begin{equation*}
P\om(X,Y)
=2\Big[\int_0^1\frac{d}{dt}\ka(X,d_xY)dt\Big]
%=\Big[\int_0^1\Big(\ka\Big(\frac{d}{dt}X,d_xY\Big)
%+\ka\Big(X,\frac{d}{dt}d_xY\Big)\Big)dt\Big]\\
%=\Big[\int_0^1\Big(\ka\Big(\frac{d}{dt}X,d_xY\Big)
%-\ka\Big(d_xX,\frac{d}{dt}Y\Big)
%+d_x\ka\Big(X,\frac{d}{dt}Y\Big)\Big)dt\Big]
=2\Big[\int_IdX\wedge_\ka dY\Big].
\end{equation*}

A 2-cocycle on the path group $PG$ integrating $P\om$ 
can be obtained from (\ref{formula2}) using the invariance of $\ka$, 
formulas from the appendix and the fact that $\de_x^l g(0,x)=0$
for $g\in PG$ and $x\in M$ by the following computation:
\begin{align*}
C(f,g)
%&=\int_0^1\Big(\int_0^s\om(\Ad(g(t)^{-1})\de^lf(s),\de^lg(t))dt\Big)ds\\
&\stackrel{(\ref{formula2})}{=}
2\int_0^1\Big(\int_0^s[\ka(\Ad (g(t,x)^{-1})i_{\partial_s}\de^lf(s,x),
d_xi_{\partial_t}\de^lg(t,x))]dt\Big)ds\\
&\stackrel{(\ref{21})}=2\int_0^1\Big(\int_0^s[\ka(i_{\partial_s}\de^lf(s,x),
\Ad(g(t,x))(\frac{d}{dt}\de^l_xg(t,x)\\
&\qquad\qquad\qquad\qquad\qquad\qquad
+[i_{\partial_t}\de^lg(t,x)),
\de^l_xg(t,x)])]dt\Big)ds\\
&\stackrel{(\ref{ad})}{=}2\int_0^1\Big(\int_0^s\Big[\ka(i_{\partial_s}\de^lf(s,x),
\frac{d}{dt}(\Ad (g(t,x))\de^l_xg(t,x))\Big]dt\Big)ds\\
&=2\int_0^1[\ka(i_{\partial_s}\de^lf(s,x),
\Ad(g(s,x))\de^l_xg(s,x)]ds=2\int_0^1[\ka(i_{\partial_s}\de^lf,\de^r_xg)]ds.
\end{align*}

The group $\Om_0G$ of null-homotopic loops in the current group $G$
can be identified with the space of those smooth maps $g\in C^\oo(I\x M,H)$
with $g(0,\cdot)=g(1,\cdot)=e$ for which there exists 
a smooth homotopy $(s,t,x)\in I\x I\x M\mapsto\bar g(s,t,x)\in H$
with $\bar g(0,\cdot,\cdot)=e$, $\bar g(1,\cdot,\cdot)=g$
and $\bar g(\cdot,0,\cdot)=\bar g(\cdot,1,\cdot)=e$.
The homotopy class of $\bar g$ is identified with an element of $\widetilde{\Om_0G}$.
The resolving map for the cocycle $C$ defined on $\widetilde{\Om_0G}$ is:
\begin{equation*}
\La([\bar g])\stackrel{(\ref{big})}{=}-\int_{\bar g}\om^{l}
=-2\Big[\int_0^1\int_0^1\ka(i_{\partial_s}\de^l\bar g,
di_{\partial_t}\de^l\bar g)dsdt\Big]
\end{equation*}

These are not yet the type of expressions expected from \cite{M89} and \cite{LMNS98}.
To get them, we will apply the symmetrization procedure from Remark \ref{prime}.
A 2-cocycle on the path group of the current group $G$, integrating $P\om$ is
\begin{multline*}
C_{\sym}(f,g)
%=C(f,g)+C'(f,g)
=\frac12(C(f,g)-C(g^{-1},f^{-1}))\\
=\int_0^1[\ka(i_{\partial_s}\de^lf,\de_x^rg)]ds
-\int_0^1[\ka(i_{\partial_s}\de^rg,\de^l_xf)]ds=\Big[\int_I\de^lf\wedge_\ka\de^rg\Big],
\end{multline*}
obtained by fiber integrating the 2-form 
$\de^lf\wedge_\ka\de^rg\in\Om^2(I\x M,V)$.
At the last step we used the identity  $\de^lf\wedge_\ka\de^rg
=dt\wedge(\ka(i_{\partial_t}\de^lf,\de^r_xg)-\ka(i_{\partial_t}\de^rg,\de^l_xf))
+\de^l_xf\wedge_\ka\de^r_xg$ obtained from the relation 
$\de^lf=(i_{\partial_t}\de^lf)dt+\de^l_xf
\in\Om^1(M\x I,\h)$.

We define $\et=\frac{1}{6}\th^l\wedge_\ka(\th^l\wedge_{[,]}\th^l)
%=\frac1{12}\ka(\th^l,[\th^l,\th^l])
\in\Om^3(H,V)$, where $\th^l\in\Om^1(H,\h)$ is the left Maurer-Cartan form.
Like the Cartan 3-form on a simple Lie group, $\et$ is a closed biinvariant 3-form
and $\et(X,Y,Z)=\ka(X,[Y,Z])$ for $X,Y,Z\in\h$.
The resolving map for the cocycle $C_{\sym}$ is
\begin{align*}
\La_{\sym}([\bar g])
%&=\frac12(\La([\bar g])-\La([\bar g^{-1}]))\\
&\stackrel{(\ref{bigprime})}{=}\frac12\La([\bar g])
+\Big[\int_0^1\int_0^1\ka(\Ad(\bar g)i_{\partial_s}\de^l\bar g,d\Ad(\bar g)
i_{\partial_t}\de^l\bar g)dsdt\Big]\\
&\stackrel{(\ref{ad})}{=}
%\frac12\La(\bar g)+\Big[\int_0^1\int_0^1
%\ka(\Ad(\bar g)i_{\partial_s}\de^l\bar g,
%\Ad(\bar g)(di_{\partial_t}\de^l\bar g
%+[\de^l\bar g,i_{\partial_t}\de^l\bar g]))dsdt\Big]
%\\&=
-\Big[\int_0^1\int_0^1\ka(i_{\partial_s}\de^l\bar g,[i_{\partial_t}\de^l\bar g,
\de^l\bar g])dsdt\Big]\\
&=-\frac1{6}\Big[\int_{I\x I}\de^l\bar g\wedge_\ka(\de^l\bar g\wedge_{[,]}
\de^l\bar g)\Big]
=-\Big[\int_{I\x I}\bar g^*\et\Big],
\end{align*} 
for the homotopy $\bar g:I\x I\x M\to H$.

This gives a constructive proof for the following slight generalization 
of the result from \cite{LMNS98} mentioned in Remark \ref{discrete}.

\begin{theo}
Assuming the period group $\Pi_\om\subset\a$ of the continuous Lie  algebra
2-cocycle (\ref{mn}) is discrete,
let $\Ga\subset\a$ be a discrete set containing $\Pi_\om$ and $\exp:\a\to A=\a/\Ga$.
%Applying Corollary \ref{second}
The $A$-valued group cocycle $c_{\sym}=\exp\o C_{\sym}$ on $PG$,
$$
c_{\sym}(f,g)=\exp\Big[\int_I\de^lf\wedge_\ka\de^rg\Big],
$$
%with $C,C'$ given by (\ref{formula}), (\ref{cprime})
integrates the Lie algebra cocycle $P\om$. 
The restriction of $c_{\sym}$ to $\Om_0G$ is the coboundary of
the inverse of $\la_{\sym}=\exp\o\La_{\sym}:\Om_0G\to A$,
$$
\la_{\sym}(g)=\Big(\exp\Big[\int_{I\x I}\bar g^*\et\Big]\Big)^{-1}.
$$
%with $\La,\La'$ given by (\ref{big}) and (\ref{bigprime}), 
The quotient group 
$(A\times_{c_{\sym}} PG)/\Graph(\la_{\sym})$ is an abelian Lie group extension 
of the universal covering group $\tilde G$ of the current group $G$
by $A$, integrating $\om$.
\end{theo}

%There exist also central extensions of $G$ itself
%integrating $\om$ \cite{MN03}.

\begin{rema}
We apply this theorem to the special case treated in \cite{LMNS98},
where $\ka$ is a multiple of the Killing form of the simple Lie algebra $\h$,
so $\et$ is a multiple of the Cartan 3-form, the constant factor being chosen
such that $\et$ is integral.
From the relation $\Pi_{\om}=\La_{\sym}(\pi_1(\Om_0G))$
in Remark \ref{periods} applied to 
$\La_{\sym}([\bar g])=-\Big[\int_{I\x I}\bar g^*\et\Big]$
follows that the period group 
$\Pi_\om=Z^1_\ZZ(M)/dC^\oo(M)$, as mentioned in Remark \ref{discrete}.

Denoting by $m:H\x H\to H$ the group multiplication and by $\pr_1,\pr_2:H\x H\to H$
the canonical projections, the Polyakov-Wiegmann formula is
$$
m^*\et=\pr_1^*\et+\pr_2^*\et-d(\pr_1^*\th^l\wedge_\ka\pr_2^*\th^r).
$$
Integrating the pullback 
of the Polyakov-Wiegmann formula by the map $(\bar f,\bar g):I\x I\x M\to H\x H$ over $I\x I$,
provides another proof that $\La_{\sym}$ resolves $C_{\sym}$.
\end{rema}

%%%%%%%%%%%%%%%%%%%%

\section{Coupled cocycle}\label{coupled}

The coupled cocycle on the current Lie algebra $\g=C^\oo(M,\h)$
was defined in \cite{N07}.
It is built with a continuous invariant symmetric bilinear form
$\ka:\h\x\h\to V$ whose image under the Cartan map
$$
\Ga:S^2(\h,V)^\h\to Z^3(\h,V),\quad\Ga(\ka)(X,Y,Z)=\ka([X,Y],Z)
$$
is a coboundary, \ie there is a 2-cochain $\be\in C^2(\h,V)$ such that
\begin{equation}\label{star}
\Ga(\ka)=d_\h\be.
\end{equation}
The corresponding coupled cocycle is 
\begin{gather*}
\om:\g\x\g\to\Om^1(M,V),\quad\om(X,Y)=\ka(X,dY)-\ka(Y,dX)-d(\be(X,Y)),
\end{gather*}
a Lie algebra 2-cocycle on $\g$ with values in the trivial module $\Om^1(M,V)$.

In \cite{N07} it is shown that the period map
$\per_\om:\pi_2(G)\to\Om^1(M,V)$ of the coupled cocycle vanishes for $G=C^\oo(M,H)_0$.
Hence there exists a central extension of $\tilde G$ by $\Om^1(M,V)$ integrating $\om$.
We explicitly realize this Lie group extension with the path group method.

\begin{rema}
The coupled cocycle is a lift to $\Om^1(M,V)$ of the cocycle
$$
\g\x\g\to\Om^1(M,V)/C^\oo(M,V)=\a,\quad
(X,Y)\mapsto 2[\ka(X,dY)]
$$ 
studied in the previous section.
In this special case when $\Ga(\ka)$ is exact, \ie a coboundary in $B^3(\h,V)$, the period map of this $\a$-valued cocycle also vanishes.
\end{rema}

The pullback cocycle $P\om$ on $P\g$ of the coupled cocycle $\om$ 
integrates to a group cocycle on the path group $PG$ because the path group is contractible. 
%An expression for such a cocycle can be obtained from (\ref{formula2}), 
%using the invariance of $\ka$ and formulas from the appendix.
A computation starting from (\ref{formula2}), similar to the one 
in the previous section, together with the symmetrization procedure,
gives the $\Om^1(M,V)$-valued cocycle $C_{\sym}$ on $PG$ as the sum of two cocycles,
one has values in the subspace of exact $V$-valued 1-forms on $M$,
the other one is $(f,g)\mapsto \int_I\de^lf\wedge_\ka\de^rg$.

More precisely, for $f,g\in PG$ we define
\begin{equation}\label{a}
m(s, t)=\ka(\de_t^rg(t,\cdot),\de_s^lf(s,\cdot))\in C^\oo(M,V)
\end{equation}
and
\begin{equation}\label{b}
C_\be(f,g)=\int_0^1\Big(\int_0^s\be(\Ad(g(t,\cdot)^{-1})\de_s^lf(s,\cdot),
\de_t^lg(t,\cdot))dt\Big)ds.
\end{equation} 
$C_\be$ is not a group cocycle in general, nevertheless we consider its 
symmetrized version $C_{\be,\sym}$ as in (\ref{csym}). 
Then, by a computation which can be found in the appendix,
\begin{equation}\label{ccc}
C(f,g)
=2\int_0^1\ka(\de_x^rg(s,\cdot),\de_s^lf(s,\cdot))ds
-d_x\Big(C_\be(f,g)+\int_0^1\int_0^s m(s, t)dtds\Big),
\end{equation}
so
\begin{gather*}
C_{\sym}(f,g)
=\int_I\de^lf\wedge_\ka\de^rg-d_x\Big(C_{\be,\sym}(f,g)
+\frac12\int_0^1\int_0^s (m(s, t)-m(t,s))dtds\Big).
\end{gather*}

Let $g:I\x M\to H$ in $\Om_0G$ and $\bar g\in C^\oo(I\x I\x M,H)$
any homotopy with $\bar g(0,\cdot,\cdot)=e$, $\bar g(1,\cdot,\cdot)=g$
and $\bar g(\cdot,0,\cdot)=\bar g(\cdot,1,\cdot)=e$.
The computation of the associated resolving map $\La$ presented in the appendix
gives
\begin{equation}\label{resolv}
\La([\bar g])=\int_0^1\ka(\de^l_xg,\de^l_tg)dt
+\int_0^1\be(\de^l_xg,\de^l_tg)dt.
\end{equation}
It follows that
$$
\La_{\sym}([\bar g])=\frac12\int_0^1\be(\de^l_xg,\de^l_tg)dt
-\frac12\int_0^1\be(\de^r_xg,\de^r_tg)dt=\La_{\sym}(g),
$$
is a map depending only on the endpoint $g$ of the homotopy class $[\bar g]$. Hence 
the map $\La_{\sym}$ descends to $\Om_0G$ and $\Pi_{\om}=\La_{\sym}(\pi_1(\Om_0G))=0$,
this giving another proof that the period group of the coupled cocycle $\om$ vanishes.
The restriction of $C_{\sym}$ to the subgroup $\Om_0G$
of null-homotopic loops based at $e$ 
is the coboundary of $-\La_{\sym}:\Om_0G\to\Om^1(M,V)$.

\begin{theo}
The quotient group $(\Om^1(M,V)\x_{C_{\sym}}PG)/\Graph(\La_{\sym})$ 
is a central extension of the universal cover $\tilde G$
of the current group
by $\Om^1(M,V)$, integrating the coupled cocycle $\om$.
Here the group cocycle $C_{\sym}$ on $PG$ is
$$
C_{\sym}(f,g)
=\int_I\de^lf\wedge_\ka\de^rg-d_x\Big(C_{\be,\sym}(f,g)
+\frac12\int_0^1\int_0^s (m(s, t)-m(t,s))dtds\Big)
$$
for $m$ and $C_\be$ given by (\ref{a}) and (\ref{b}),
and the smooth map $\La_{\sym}:\Om_0G\to\Om^1(M,V)$ is
$$
\La_{\sym}(g)=\frac12\int_0^1\be(\de^l_xg,\de^l_tg)dt
-\frac12\int_0^1\be(\de^r_xg,\de^r_tg)dt.
$$
\end{theo}

%%%%%%%%%%%%%%%%%%%%%%%%%%%%%%%%%%%%%%%%%%%%

\section{The group of volume preserving diffeomorphisms}

On the compact manifold $M$ we consider an integral volume form $\mu$.
Let $G=\Diff_\mu(M)_0$ be the connected component of the group
of volume preserving diffeomorphisms
and $\g=\X_\mu(M)$ its Lie algebra, the Lie algebra 
of divergence free vector fields \cite{KM97} Section 43.
Its subgroup, the group of exact volume preserving diffeomorphisms, 
is a Lie group 
with Lie algebra $\{X\in\X(M):i_X\mu \text{ is an exact differential form}\}$,
kernel of a flux homomorphism \cite{B97}.

Given $\et$ a closed integral 2-form on $M$,
the Lichnerowicz cocycle 
$$
\om:\g\x\g\to\RR,\quad\om(X,Y)=\int_M\et(X,Y)\mu
$$
is a Lie algebra 2-cocycle on the Lie algebra of divergence free vector fields. Indeed, by the closedness of $\et$,
$$
\sum_{cycl}\int_M\et([X_1,X_2],X_3)\mu
=\sum_{cycl}L_{X_1}(\et(X_2,X_3))\mu=0.
$$

The Lichnerowicz cocycle integrates to a central Lie group extension of 
the subgroup of exact volume preserving diffeomorphisms
\cite{I96} \cite{HV04}. It integrates also to a central Lie group extension
of the universal cover $\tilde G$ of the group of volume preserving diffeomorphisms:
the existence is proved in \cite{N04} and a construction with Kostant's prequantization
central extension is given in \cite{V05}.
In this section we use the method of Section 3 to explicitly realize this
central Lie group extension.

Because $\et$ is integral, there exists a principal circle bundle 
$q:(P,\th)\to(M,\et)$ with connection 1-form $\th$
and curvature 2-form $\et$, so that $q^*\et=d\th$.
There is a natural volume form on $P$ determined by $\mu$ and $\th$,
namely $\tilde\mu=\th\wedge q^*\mu$. 
It has the property that for any $f\in C^\oo(M)$,
$\int_P(q^*f)\tilde\mu=\int_Mf\mu$.
Cosidering the principal $\TT$-action on $P$,
each $\TT$-invariant divergence free vector field 
in $\X_{\tilde\mu}(P)^\TT$ projects to a divergence free vector field in $\X_\mu(M)$.
Every vector field $X\in\X(M)$ has a horizontal lift to $P$, denoted
by $X^{hor}$, uniquely defined by $\th(X^{hor})=0$ and $q_*X^{hor}=X$.
Moreover, if $X$ is divergence free \wrt $\mu$, then $X^{hor}$ is divergence free 
\wrt $\tilde\mu$, hence the horizontal lift provides a section of
the abelian Lie algebra extension
\begin{equation}\label{gauge}
0\to C^\oo(M)\to\X_{\tilde\mu}(P)^\TT\to\X_\mu(M)\to 0.
\end{equation}

Let $g$ be a volume preserving diffeotopy 
of $M$ starting at the identity. One can lift it to the
volume preserving diffeotopy $g^{hor}$ of $P$, 
starting at the identity and defined by 
\begin{equation}\label{111}
\de^r (g^{hor})=(\de^rg)^{hor}.
\end{equation}
%Obviously $q\o g^{hor}=g\o q$ and
The diffeotopy $g^{hor}$ consists of $\TT$-equivariant diffeomorphisms of $P$.
For $X\in\X_\mu(M)$ and $t\in I$ the vector field $g^{hor}(t)^*X^{hor}$ on $P$
descends to the vector field $g(t)^*X$, but is not necesarily horizontal.
The failure of horizontality is measured by 
$$
g^{hor}(t)^*X^{hor}-(g(t)^*X)^{hor}=q_*(\th(g^{hor}(t)^*X^{hor}))\in C^\oo(M).
$$
Here the $\TT$-equivariance of $g^{hor}(t)$ together with the
$\TT$-invariance of $X^{hor}$ assure that
$\th(g^{hor}(s)^*X^{hor})$ is the pullback of a function on $M$,
function denoted by $q_*(\th(g^{hor}(t)^*X^{hor}))$.
The computation of a group 2-cocycle $C$ on $PG$ 
integrating $P\om$, using the formula (\ref{formula2}), 
gives $C$ as the integral over $M$ of an expression of this type.

Indeed, using the fact that
the adjoint action in $\Diff(M)$ is $\Ad(g)X=(g^{-1})^*X=Tg\o X\o g^{-1}$ in the first step,
the relation $\de^l g=\Ad(g^{-1})\de^rg=g^*\de^rg$ following from
(\ref{100}) in the second step, and the fact that the horizontal lift of 
a divergence free vector field is again divergence free in step 4, we get.
\begin{align*}
C(f,g)
%&=\int_0^1\int_0^s\om(\Ad(g(t)^{-1})\de^lf(s),\de^lg(t))dtds\\
&\stackrel{(\ref{formula2})}{=}
\int_0^1\int_0^s\Big(\int_M\et(g(t)^*\de^l f(s),\de^l g(t))\mu\Big) dtds\\
&=\int_0^1\int_0^s\Big(\int_M (g(t)^{-1*}\et)(\de^l f(s),\de^r g(t))\mu\Big) dtds\\
%&=\int_0^1\int_0^s\Big(\int_P\Big( q^*(g(t)^{-1*}\et)
%(\de^l f(s),\de^r g(t))\Big)\tilde\mu\Big) dtds\\
&=\int_0^1\int_0^s\Big(\int_P (dg^{hor}(t)^{-1*}\th)(\de^l f(s)^{hor},
\de^r g(t)^{hor})\tilde\mu\Big) dtds\\
%&=-\int_0^1\int_0^s\Big(\int_P (g^{hor}(t)^{-1*}\th)([\de^l f(s)^{hor},
%\de^r g(t)^{hor}])\tilde\mu\Big) dtds\\
&=\int_0^1\int_0^s\Big(\int_P \th(g^{hor}(t)^*L_{\de^r g(t)^{hor}}
\de^l f(s)^{hor})\tilde\mu\Big) dtds\\
&\stackrel{(\ref{111})}{=}
\int_0^1\int_0^s\Big(\int_P \th(\frac{d}{dt}(g^{hor}(t)^*\de^l f(s)^{hor}))
\tilde\mu\Big) dtds\\
%&=\int_0^1\Big(\int_P \th(g^{hor}(s)^*\de^l f(s)^{hor})\tilde\mu\Big) ds\\
&=\int_0^1\Big(\int_M q_*(\th(g^{hor}(s)^*\de^l f(s)^{hor}))\mu\Big) ds.
\end{align*}

Let $\bar g\in\widetilde{\Om_0G}$ be the homotopy class of a path in 
$\Om_0G$ viewed as a map $\bar g:I\x I\to G$ and 
$\hat g:I\x I\x M\to M$,
$\hat g(s,t,x)=\bar g(s,t)(x)$.
Let $p:I\x I\x M\to M$ denote the projection on the third factor.
The resolving map for $C$ is 
\begin{align*}
\La([\bar g])&\stackrel{(\ref{big})}{=}-\int_{\bar g}\om^{eq}
%=-\int_0^1\int_0^1\om(i_{\partial_s}\de^l\bar g,i_{\partial_t}\de^l\bar g)dsdt\\&
=-\int_0^1\int_0^1\Big(\int_M\et(i_{\partial_s}\de^l\bar g,
i_{\partial_t}\de^l\bar g)\mu\Big) dsdt\\
&=-\int_0^1\int_0^1\Big(\int_M\et\wedge i_{{\partial_t}\de^l\bar g}
i_{{\partial_s}\de^l\bar g}\mu\Big) dsdt=-\int_{I\x I\x M}p^*\et\wedge\hat g^*\mu.
\end{align*}

\begin{rema}
The restriction of $\La$ to $\pi_2(G)=\pi_1(\Om_0G)\subset \widetilde{\Om_0G}$
is minus the period map of $\om$ by Remark \ref{periods}. 
One can now easily see that the period group of $\om$ is discrete.
Each piecewise smooth representative 
$\ta:S^2\to G$ of a class $[\ta]\in\pi_2(G)$ defines a map 
$\hat\ta:S^2\x M\to M$.
Then $\per_\om([\ta])=-\La([\ta])=\int_{S^2\x M}p^*\et\wedge\hat\ta^*\mu\in\ZZ$,
because both $\et$ and $\mu$ are integral forms on $M$.
Hence the period group $\Pi_\om$ is contained in $\ZZ$. 
\end{rema}

Let $\exp:\RR\to \TT\cong \RR/\ZZ$,
$c=\exp\o C$ and $\la=\exp\o\La$.
The construction of Theorem \ref{main} yields a central extension
integrating the Lichnerowicz cocycle:
\comment
{
\begin{theo}
Let $\mu$ be an integral volume form on a compact manifold $M$
and $\et$ a closed integral 2-form.
The Lichnerowicz cocycle $\om(X,Y)=\int_M\et(X,Y)\mu$
on $\X_\mu(M)$ integrates to the 1-dimensional 
central extension of the universal covering group $\tilde G$ of $G=\Diff_\mu(M)_0$,
obtained as the quotient group $(\TT\x_c PG)/\Graph(\la)$.
The group 2-cocycle $c:PG\x PG\to \TT$ is
$$
c(f,g)=\exp\int_0^1\Big(\int_Mq_*(\th(g^{hor}(t)^*\de^l f(t)^{hor}))\mu\Big) dt,
$$
for $\th$ a connection 1-form on a principal circle bundle $q:P\to M$ 
with curvature $\et$, $X^{hor}$ the horizontal lift of $X$, 
$g^{hor}$ the diffeotopy of $P$ starting at the identity and defined by
$\de^r (g^{hor})=(\de^r g)^{hor}$.
The map $\la:\Om_0G\to \TT$ is
$$
\la(g)=\Big(\exp\int_{I\x I\x M}p^*\et\wedge\hat g^*\mu\Big)^{-1},
$$
for $\bar g$ a path from $\id_M$ to $g$ in $\Om_0G$ defining $\hat g:I\x I\x M\to M$,
$\hat g(s,t,x)=\bar g(s,t)(x)$
and $p:I\x I\x M\to M$ the projection on the third factor.
\end{theo}
}

\begin{theo}
The quotient group $(\TT\x_c PG)/\Graph(\la)$ 
is a central extension 
of the universal covering group $\tilde G$ of $G=\Diff_\mu(M)_0$
by the circle $\TT$, integrating the Lichnerowicz cocycle $\om(X,Y)=\int_M\et(X,Y)\mu$
defined with the closed integral 2-form $\et$ on $M$.
Here the group cocycle $c:PG\x PG\to \TT$ is
$$
c(f,g)=\exp\int_0^1\Big(\int_Mq_*(\th(g^{hor}(t)^*\de^l f(t)^{hor}))\mu\Big) dt,
$$
and the smooth map $\la:\Om_0G\to \TT$ is
$$
\la(g)=\Big(\exp\int_{I\x I\x M}p^*\et\wedge\hat g^*\mu\Big)^{-1}.
$$ 
\end{theo}

If $M$ is 2-dimensional, then $\mu$ can be viewed as a symplectic form and $\g=\X_\mu(M)$ as the Lie algebra of symplectic vector fields. The kernel of the infinitesimal flux homomorphism $\flux_\mu:X\in\X_\mu(M)\mapsto[i_X\mu]\in H^1_{dR}(M)$ is the Lie algebra of Hamiltonian vector fields. The Lie algebra cohomology class
$[\om]\in H^2(\g)$ 
of the Lichnerowicz cocycle is the pullback by $\flux_\mu$ of a multiple of the skew-symmetric pairing 
$(a,b)\mapsto\int_M a\wedge b$ on $H^1_{dR}(M)$. This is a consequence of the fact that $\om$ is cohomologous to a multiple of the 2-cocycle 
$(X,Y)\mapsto\int_M\mu(X,Y)\mu$ on $\g$, whose restriction to the Lie algebra of Hamiltonian vector fields is trivial. 

For an arbitrary compact symplectic manifold $M$,
Lie algebra 2-cocycles on the Lie algebra of Hamiltonian vector fields, having non-zero cohomology classes, are associated to closed 1-forms on $M$ \cite{R95}. It seems that this path method does not work for them.

%%%%%%%%%%%%%%%%%%%%%%%

\renewcommand{\theequation}{A\arabic{equation}}
\setcounter{equation}{0}

%%%%%%%%%%%%%%%%%%%%%

\section*{Appendix: Logarithmic derivatives}

Let $G$ be a Lie group with Lie algebra $\g$.
The left logarithmic derivative of a smooth map $g:M\to G$
is the 1-form $\de^lg\in\Om^1(M,\g)$, $\de^lg(X_x)=g(x)^{-1}T_xg.X_x$
for any $X_x\in T_xM$. It is the pull-back of the left Maurer-Cartan form
$\th^l\in\Om^1(G,\g)$ by the map $g$. 
When $M=I=[0,1]$, we identify
the left logarithic derivative of a smooth curve $g:I\to G$
with the curve $i_{\partial_t}\de^lg(t)=g(t)^{-1}\dot g(t)$,
usually also denoted by $\de^lg:I\to\g$.

The right logarithmic derivative is defined similarly and
\begin{gather}\label{100}
\de^rg=\Ad(g)\de^lg=-\de^l(g^{-1}).
\end{gather}

For $f,g:M\to G$, $\xi:M\to\g$ smooth maps and
$\Ph:G\to H$ a Lie group homomorphism with derivative $\ph:\g\to\h$,
the following formulae hold:
\begin{gather}
%\de^r(g^{-1})=-\de^lg\\
%\de^rg=\Ad(g)\de^lg\\
\de^l(fg)=\de^l(g)+\Ad(g)^{-1}\de^l(f)\label{prod}\\
d\de^lg+\frac12\de^lg\wedge_{[,]}\de^lg=0\label{mc}\\
d\Ad(g)\xi=\Ad(g)[d\xi+\ad(\de^lg)\xi]\label{ad}
%=\Ad(g)d\xi+\ad(\de^rg)(\Ad(g)\xi)\\
\\
\de^l(\Ph\o g)=\ph\o\de^lg,
\end{gather}

The right Maurer-Cartan equation (\ref{mc}) applied to $X,Y\in\X(M)$
becomes 
\begin{equation}\label{300}
(d\de^lg)(X,Y)=[\de^lg(Y),\de^lg(X)].
\end{equation}
For $g:I\x M\to G$, $X\in\X(M)$
and $Y=\partial_t$ it implies
\begin{equation}\label{21}
%\frac{d}{dt}i_{\partial_s}\de^lf-\frac{d}{ds}i_{\partial_t}\de^lf=
%[i_{\partial_s}\de^lf,i_{\partial_t}\de^lf]\\
d_x(\de_t^lg)=\frac{d}{dt}\de^l_xg+[\de_t^lg,\de^l_xg].
\end{equation}
Here $\de_x^l$ and $d_x$ denote the logarithmic derivative and the differential
on $M$ considering $t\in I$ as a parameter, and $\de_t^lg=i_{\partial_t}\de^lg$.
 
From (\ref{ad}) and (\ref{21}) follows that
\begin{equation}\label{199}
\frac{d}{dt}\de^r_xg=\Ad(g)d_x\de^l_tg.
\end{equation}

The right logarithmic derivative satisfies the left Maurer-Cartan equation, so
the analogue of (\ref{21}) for the right logarithmic derivative is
\begin{equation}\label{200}
d_x(\de_t^rg)=\frac{d}{dt}\de^r_xg-[\de_t^rg,\de^r_xg]
\end{equation}
and from (\ref{ad}) follows
\begin{equation}\label{101}
d_x\Ad(g^{-1})\xi=\Ad(g^{-1})d_x\xi-[\de^l_xg,\Ad(g^{-1})\xi].
\end{equation}

%%%%%%%%%%%%%%%%%%%%%%%%%%%%%%%%

\end{document}